\input amstex 
\documentstyle{amsppt} 
\loadbold
\magnification=1200
\NoBlackBoxes

\pagewidth{6.43truein}
\pageheight{8.5truein}

\define\A{{\Cal A}}
\define\C{{\Bbb C}}
\redefine\D{{\Bbb D}}

\define\R{{\Bbb R}}
\define\T{{\Bbb T}}
\redefine\phi{{\varphi}}
\redefine\epsilon{{\varepsilon}}
\define\ac{\acuteaccent}

\define\pr{\operatorname{pr}}
\define\set#1{\{#1\}}
\redefine\cdot{{\boldsymbol\cdot}}

\hyphenation{pluri-sub-har-mon-ic}

\refstyle{B}
\NoRunningHeads
\TagsOnRight

\topmatter 
\title Plurisubharmonicity of envelopes  \\ 
of disc functionals on manifolds \endtitle
\author Finnur L\ac arusson and Ragnar Sigurdsson \endauthor
\address Department of Mathematics, University of Western Ontario,
London, Ontario N6A~5B7, Canada \endaddress
\email larusson\@uwo.ca \endemail
\address Science Institute, University of Iceland, Dunhaga 3,
IS-107 Reykjav\ac ik, Iceland \endaddress
\email ragnar\@hi.is \endemail 

\thanks The first-named author was supported in part by the Natural
Sciences and Engineering Research Council of Canada.  \endthanks

\thanks Completed 19 November 2001; minor changes 24 January 2002.
\endthanks

\abstract We show that a disc functional on a complex manifold has a
plurisubharmonic envelope if all its pullbacks by holomorphic
submersions from domains of holomorphy in affine space do and it is
locally bounded above and upper semicontinuous in a certain weak sense. 
For naturally defined classes of disc functionals on manifolds, this
result reduces a property somewhat stronger than having a
plurisubharmonic envelope to the affine case.  The proof uses a recent
Stein neighbourhood construction of Rosay, who proved the
plurisubharmonicity of the Poisson envelope on all manifolds.  As a
consequence, the Riesz envelope and the Lelong envelope are
plurisubharmonic on all manifolds; for the former, we make use of new
work of Edigarian.  The basic theory of the three main classes of disc
functionals is thereby extended to all manifolds.  \endabstract

\subjclass Primary: 32U05; secondary: 32Q99, 32U35  \endsubjclass

\endtopmatter

\document

\specialhead Introduction \endspecialhead

\noindent The theory of disc functionals was founded by Evgeny A\. 
Poletsky in the late 1980s with the papers \cite{P1, PS}.  In its
applications in pluripotential theory and elsewhere in complex analysis,
we rarely consider a single disc functional in isolation.  Rather, we
normally have a class of disc functionals of a given type on each
complex manifold, preserved by taking pullbacks by holomorphic maps. 
The question of plurisubharmonicity of envelopes is fundamental.  An
analytic disc can always be lifted by a suitable holomorphic submersion
from a domain of holomorphy in affine space, and, with varying degrees of
effort, so can certain configurations of discs.  This suggests the
possibility of reducing plurisubharmonicity of envelopes of disc
functionals of a given type to the affine case, where the problem may
have been studied already. 
 
The main theorem of this paper states that a property somewhat stronger
than having a plurisubharmonic envelope reduces to the affine case,
meaning that a disc functional has this property if all its pullbacks by
holomorphic submersions from domains of holomorphy in affine space do. 
We do not know if the property of having a plurisubharmonic envelope by
itself reduces to the affine case.  Still, as a consequence,
plurisubharmonicity of envelopes is established for the three main
classes of disc functionals, the Poisson, Riesz, and Lelong functionals,
on all manifolds. 

Thereby, an extra hypothesis is removed from various theorems based on
the theory of disc functionals, such as the product property of
relative extremal functions \cite{EP}, the product property of
generalized Green functions \cite{LS2}, and results in Chapters 6 and 7
of \cite{LS1}.  For instance, our Kontinuit\"atssatz, characterizing
pseudoconvexity in terms of analytic discs, is valid for all manifolds. 

In our paper \cite{LS1}, we identified the construction of a certain
Stein neighbourhood as the obstacle to extending the theory of disc
functionals to all manifolds.  The long-awaited breakthrough came with
Rosay's paper \cite{Rs} in April of 2001, in which the Poisson envelope
was shown to be plurisubharmonic on all manifolds.  We thank Armen
Edigarian for sharing with us a draft of a simplified version of Rosay's
argument \cite{E2}, on which the proof of Theorem 1.2 is partly based.

For background information on the Poisson, Riesz, and Lelong functionals
and their envelopes, we refer the reader to \cite{LS1, LS2}. 

Some notation: We denote by $D_r$ the open disc
$\set{z\in\C\,;\,|z|<r}$, $r>0$, by $\overline\R$ the extended real line
$[-\infty,+\infty]$, and by $\lambda$ the {\it normalized} arc length
measure on the unit circle $\T$ (note that in our previous papers, we
did not normalize $\lambda$).

\specialhead 1. Reduction to the affine case  \endspecialhead

\noindent Before stating and proving our reduction theorem, we recall a
few definitions.  Let $X$ be a complex manifold and $\A_X$ be the set of
maps $f:\overline\D\to X$ which are holomorphic in a neighbourhood of
the closure $\overline\D$ of the unit disc $\D$.  Such maps are called
{\it (closed) analytic discs\,} in $X$.  For a point $p\in X$, we denote
the constant disc at $p$ by $p$.  A {\it disc functional\,} on $X$ is a
map $H:\Cal A_X\to\overline\R$.  The {\it envelope\,} of $H$ is the
function $EH:X\to\overline\R$ defined by the formula
$$
EH(x)=\inf\,\set{H(f)\,;\, f\in\A_X, f(0)=x}, \qquad x\in X. 
$$ 
If $\phi:Y\to X$ is a holomorphic map, then the pullback $\phi^*H$ is
the disc functional on $Y$ defined by the formula $f\mapsto H(\phi\circ
f)$.  Clearly, $EH\circ\phi\leq E\phi^*H$.  

\proclaim{1.1.  Lemma} Let $H$ be a disc functional on an
$n$-dimensional complex manifold $X$ such that the envelope $E\phi^*H$
is upper semicontinuous for every holomorphic submersion $\phi$ from the
$(n+1)$-dimensional polydisc into $X$.  Then $EH$ is upper semicontinuous.
\endproclaim

\demo{Proof} We first observe that $EH$ is nowhere $+\infty$.  Namely,
for $p\in X$, choose a coordinate polydisc $U$ centred at $p$ and let
$\phi$ be the projection $U\times\D\to U$.  Since $E\phi^*H$ is upper
semicontinuous by hypothesis, $E\phi^*H(p,0)<+\infty$, so there is an
analytic disc $f$ in $U\times\D$ with $f(0)=(p,0)$ and
$\phi^*H(f)<+\infty$.  Then $\phi\circ f$ is an analytic disc in
$U\subset X$ with $0\mapsto p$, and $EH(p)\leq H(\phi\circ f)=\phi^*H(f)
<+\infty$. 

Now let $p\in X$ and take $a>EH(p)$.  We need to find a neighbourhood
$U$ of $p$ such that $EH<a$ on $U$.  By definition of $EH$, there exists
an analytic disc $f_0:D_r\to X$, $r>1$, such that $f_0(0)=p$ and
$H(f_0)<a$.  In the proof of Lemma 2.3 in \cite{LS1}, using the theorem
of Siu on the existence of Stein neighbourhoods of Stein subvarieties,
we showed that for each $t\in(1,r)$, there exists a biholomorphism
$\Psi$ from a neighbourhood of the graph $\set{(z,f_0(z))\,;\, z\in
D_t}$ in $D_t\times X$ onto $D_t^{n+1}$, such that
$\Psi(z,f_0(z))=(z,0)$ for all $z\in D_t$.  We note that this also
follows from an older result of Royden \cite{Ry}.

Consider the holomorphic submersion $\phi=\pr\circ\Psi^{-1}:D_t^{n+1}\to
X$, where $\pr:\C\times X\to X$ is the projection.  The disc $\tilde
f_0:D_t \to D_t^{n+1}$, $z\mapsto (z,0)$, is a lifting of $f_0$ by
$\phi$, so $E\phi^*H(0)\leq\phi^*H(\tilde f_0)=H(f_0)<a$.  By
assumption, there exists a neighbourhood $W$ of $0$ in $D_t^{n+1}$ such
that $E\phi^*H<a$ on $W$.  Since $\phi$ is a submersion, $U=\phi(W)$ is
an open neighbourhood of $p$ in $X$.  If $q\in U$, then there is $\tilde
q\in W$ and $\tilde f\in\A_V$ with $\phi(\tilde q)=q$, $\tilde
f(0)=\tilde q$, and $\phi^* H(\tilde f)<a$.  Then $f=\phi\circ\tilde
f\in\A_X$, $f(0)=q$, and $EH(q)\leq H(f)=\phi^*H(\tilde f)<a$. 
\qed\enddemo

The following theorem is the main result of the paper.

\proclaim{1.2.  Reduction Theorem} A disc functional $H$ on a complex
manifold $X$ has a plurisubharmonic envelope if it satisfies the
following three conditions. 
\roster
\item The envelope $E\phi^*H$ is plurisubharmonic for every holomorphic 
submersion $\phi$ from a domain of holomorphy in affine space into $X$.
\item There is an open cover of $X$ by subsets $U$ with a pluripolar subset
$Z\subset U$ such that for every $h\in\A_U$ with
$h(\overline\D)\not\subset Z$, the function $w\mapsto H(h(w))$ is
dominated by an integrable function on $\T$. 
\item If $h\in \A_X$, $w\in\T$, and $\epsilon>0$, then $w$ has a
neighbourhood $U$ in $\C$ such that for every sufficiently small closed 
arc $J$ in $\T$ containing $w$ there is a holomorphic
map $F:D_r\times U\to X$, $r>1$, such that $F(0,\cdot)=h|U$ and
$$
\frac 1 {\lambda(J)}\underline {\int_J} H(F(\cdot,t))\,d\lambda(t) \leq
EH(h(w))+\epsilon.
$$
\endroster 
\endproclaim

We view \therosteritem3 as a weak upper semicontinuity condition on $H$. 
Such conditions, subtle and complicated as they are, appear naturally 
and inevitably in the theory of disc functionals.
In the absence of any measurability assumptions on $H$, the integral on
the left-hand side of the inequality is a lower integral, i.e., the
supremum of the integrals of all integrable Borel functions dominated by
the integrand.  The following remark is trivial, but worth stating for
the record. 

\proclaim{1.3.  Remark} A disc functional $H$ on a complex manifold $X$
satisfies hypothesis \therosteritem3 above if for each $f\in\A_X$ and
$\beta>H(f)$ there exists a neighbourhood $V$ of $p=f(0)$ and a
holomorphic map $F:D_r\times V\to X$, $r>1$, such that $F(\cdot,p)=f$,
$F(0,x)=x$, and $H(F(\cdot,x))<\beta$ for all $x\in V$.  
\endproclaim

By Lemma 2.3 in \cite{LS1}, the Poisson functional for an upper
semicontinuous function and the Lelong functional satisfy the condition
in the remark, but we need the weaker condition \therosteritem3 to
capture the Poisson functional for a plurisuperharmonic function and
thereby the Riesz functional. 

We view \therosteritem2 as a mild and easily-verifiable boundedness
condition.  Note that \therosteritem2 holds if $H$ is bounded above on
the set of analytic discs with image in $K$, whenever $K$ is a compact
subset of $X$.  The Poisson functional for an upper semicontinuous
function and the Lelong functional satisfy this stronger condition, but
we need \therosteritem2 to cover the Poisson functional for a
plurisuperharmonic function.

\demo{Proof of Theorem 1.2} By \therosteritem1 and Lemma 1.1, $EH$ is upper 
semicontinuous, so it remains to show that 
$$EH(h(0)) \leq \int_\T EH\circ h \,d\lambda $$
for every $h\in\A_X$.  This follows
if for every $\epsilon>0$ and every continuous function $v:X\to\R$ with 
$v\geq EH$, there exists $g\in\A_X$ such that $g(0)=h(0)$ and
$$H(g)\leq \int_\T v\circ h \,d\lambda +\epsilon. $$
Say $h$ is holomorphic on $D_r$, $1<r<2$.  

We may assume that $h(\overline\D)$ lies in any element $U$ of an open
cover of $X$ and is not contained in a given pluripolar subset of $U$. 
Namely, suppose an upper semicontinuous function $u$ has the sub-mean
value property with respect to analytic discs that do not lie in a
pluripolar set $Z$.  Since $Z$ has measure zero, $u$ has the sub-mean
value property for balls, so $u$ is subharmonic.  By shrinking the
domain, we may assume that $Z= w^{-1}(-\infty)$, where $w<0$, $w\neq
-\infty$ is plurisubharmonic.  Then $u+cw$, $c>0$, has the sub-mean
value property with respect to all analytic discs, so it is
plurisubharmonic, and the regularization $\tilde u\leq u$ of the limit
of $u+cw$ as $c\to 0$ is plurisubharmonic and equal to $u$ off $Z$. 
Since $u$, $\tilde u$ are subharmonic and equal off a nullset, they are
equal everywhere, so $u$ is plurisubharmonic. 

By a simple compactness argument, it follows from \therosteritem2 and
\therosteritem3 that we can find finitely many closed arcs
$J_1,\dots,J_m$ in $\T$, open discs $U_j$ centred on $\T$, relatively 
compact in $D_r$, containing $J_j$ with mutually disjoint closures, and 
holomorphic maps $F_j:D_s\times U_j\to X$, $s>1$, such that 
$F_j(0,\cdot)=h|U_j$ and
$$
\underline {\int_{J_j}} H(F_j(\cdot,t))\,d\lambda(t) \leq
{\int_{J_j}} v\circ h \,d\lambda +\frac \epsilon 4 \lambda(J_j),
$$
$$
\int_{\T\setminus\cup J_j}|v\circ h|\,d\lambda<\frac \epsilon 4, 
$$
and
$$
{\overline\int}_{\T\setminus\cup J_j} H(h(w))\,d\lambda(w)<\frac
\epsilon 4.
$$
Let $U_0=D_r$ and $F_0:D_s\times U_0\to X$, $(z,w)\mapsto h(w)$.

The graph $\set{(z,w,F_j(z,w))\,;\, z\in D_s, w\in U_j}$, $0\leq j\leq
m$, is a Stein submanifold of $D_s\times U_j\times X$, being
biholomorphic to $D_s\times U_j$.  By a now-familiar argument based 
on Siu's theorem, as in the proof of Lemma 2.3 in \cite{LS1}, we 
conclude that if we slightly shrink $U_j$ and $s$, then
there is a biholomorphism $\Psi_j$ from a neighbourhood of the graph
onto $D_s\times U_j\times D_s^n$ such that
$$\Psi_j(z,w,F_j(z,w))=(z,w,0), \qquad z\in D_s, w\in U_j.$$ By again
shrinking $U_j$ and $s$, and choosing $\gamma>0$ sufficiently small, we
obtain a holomorphic map $\Phi_j:U_j\times D_s^{n+3}\to\C^4\times X$,
well defined by the formula
$$\Phi_j(\zeta)=(\zeta_1,\zeta_2,-\zeta_2,-\zeta_1)+
(0,0,\Psi_j^{-1}(\zeta_2+\gamma\zeta_4,\zeta_1+\gamma\zeta_3,\zeta')),$$
where $\zeta=(\zeta_1,\dots,\zeta_4,\zeta')$ and
we intend to add the first term on the right to the first four
components of the second term.  The map $\Phi_j$ is a biholomorphism 
from $U_j\times D_s^{n+3}$ onto its image in $\C^4\times X$, and
$$\Phi_j(w,z,0)=(w,z,0,0,F_j(z,w)), \qquad w\in U_j, z\in D_s,$$
so in particular,
$$\Phi_j(w,0)=(w,0,0,0,h(w)), \qquad w\in U_j.$$
Let
$$K_0=\set{(w,0,0,0,h(w))\,;\, w\in\overline\D}=\Phi_0(\overline\D\times 
\set{0}^{n+3})$$ 
and 
$$K_j=\set{(w,z,0,0,F_j(z,w))\,;\, w\in J_j, z\in\overline\D} = \Phi_j(J_j
\times\overline\D\times\set{0}^{n+2}), \qquad j=1,\dots,m.$$

The crux of the proof is the existence of a Stein neighbourhood $V$ of
$K=K_0\cup\dots\cup K_m$ in $\C^4\times X$.  Let us assume this for a
moment and finish proving that $EH$ has the sub-mean value property. 
Let $\tau:V\to \C^N$ be an embedding and $\sigma:W\to\tau(V)$ be a
holomorphic retraction from a Stein neighbourhood $W$ of $\tau(V)$ in
$\C^N$.  Finally, let $\phi$ be the holomorphic submersion
$\pr\circ\tau^{-1}\circ\sigma:W\to X$, where $\pr:\C^4\times X\to X$ is
the projection. 

By assumption, $\tilde u=E\phi^*H$ is plurisubharmonic on $W$, so
$$\tilde u(\tilde h(0)) \leq \int_\T\tilde u\circ\tilde h\,d\lambda,$$
where $\tilde h:D_r\to W$ is the lifting $w\mapsto \tau(w,0,0,0,h(w))$
of $h$ by $\phi$.  Hence, there is a disc $\tilde g\in\A_W$
with $\tilde g(0)=\tilde h(0)$ such that 
$$H(g)=\phi^*H(\tilde g) \leq \int_\T\tilde u\circ\tilde h\,d\lambda
+ \frac \epsilon 4,$$
where $g=\phi\circ\tilde g\in\A_X$ has $g(0)=h(0)$.  
Now if $w\in J_j$, $1\leq j\leq m$,
then $z\mapsto\tau\circ\Phi_j(w,z,0)$ is a lifting by $\phi$ of
$F_j(\cdot,w)$ with $0\mapsto\tilde h(w)$, so $\tilde u(\tilde h(w))\leq
H(F_j(\cdot,w))$.  Also, if $w\in\T\setminus\cup J_j$, then $\tilde
u(\tilde h(w)) \leq \phi^*H(\tilde h(w)) = H(h(w))$.  Hence,
$$\align
\int_\T\tilde u\circ\tilde h\,d\lambda &\leq
\sum_{j=1}^m \underline{\int_{J_j}} H(F_j(\cdot,w))\,d\lambda(w)
+{\overline\int}_{\T\setminus\cup J_j} H(h(w))\,d\lambda(w) \\
&\leq \int_{\cup J_j} v\circ h\,d\lambda + \frac \epsilon 4
+\frac\epsilon 4 \leq \int_\T v\circ h\,d\lambda + \frac{3\epsilon} 4.
\endalign $$

We now turn to the all-important Stein neighbourhood.  We shall
construct a continuous, strictly plurisubharmonic exhaustion function
$\rho:V\to[0,1)$ on a neighbourhood $V$ of $K$ in $\C^4\times X$. 

For $j=1,\dots,m$, choose open discs $U_j'$, $U_j''$
concentric with $U_j$, such that 
$$J_j \subset U_j''\Subset U_j' \Subset U_j.$$ 
Since $\Phi_0(w,0)=\Phi_j(w,0)$ for all
$w\in U_j$, we have 
$$\Phi_0(\overline U_j'\times D_\epsilon^{n+3})\Subset 
\Phi_j(U_j\times D_s^{n+3}) \tag * $$ 
for $\epsilon>0$ small enough.  Our next step is to modify $\Phi_0$ so
that $\Phi_0:D_r\times D_\epsilon^{n+3}\to \C^4\times X$ is still a
biholomorphism onto its image, with a slightly smaller $r$ and
$\epsilon>0$ such that (*) still holds, and $\Phi_0(w,0)$ is still
$(w,0,0,0,h(w))$ for $w\in D_r$, which is what we need above, and so
that the derivative of the composition $\Phi_j^{-1}\circ\Phi_0$ is close
to the identity at each point $(w,0)$ with $w$ near $\overline U_j'$.  
This ensures that
$$\frac 1 A |\pi_1|^2 \leq |\pi_1\circ\Phi_j^{-1}\circ\Phi_0|^2 \leq
A|\pi_1|^2 \qquad\text{near } \overline U_j'\times D_\epsilon^{n+3},$$
with $A>1$ close to $1$ if $\epsilon>0$ is small enough.  Here,
$\pi_1:\C^{n+4}\to\C^{n+3}$ chops off the first coordinate, and
$|\cdot|$ is the euclidean norm.  This inequality with some $A>1$ is 
immediate from the Mean Value Theorem, but we need it with $A$ close to
$1$. We choose $\epsilon\in (0,s)$ 
small enough to give $\tfrac 1 A > 1-\tfrac 1 3 (1-R)^2$, where $R<1$ is the
maximum of the radii of the discs $U_1',\dots,U_m'$. 

We can modify $\Phi_0$ in this way by precomposing it with a suitable
biholomorphism $\mu$ from one neighbourhood of $D_r\times\set{0}^{n+3}$
in $D_r\times\C^{n+3}$ onto another, which is the identity on
$D_r\times\set{0}^{n+3}$.  Differentiating the inverse of
$\Phi_j^{-1}\circ \Phi_0$ at points $(w,0)$, $w\in U_j$, gives a
holomorphic, matrix-valued map on $U_1\cup\dots\cup U_m$.  We
Runge-approximate it on a neighbourhood of $\overline U_1'
\cup\dots\cup\overline U_m'$ by a holomorphic, matrix-valued map $\nu$
on $D_r$ with $\nu_{11}=1$ and $\nu_{i1}=0$ for $i\geq 2$, and let
$\mu(x)=\nu(x_1)\cdot x$, where $x$ is viewed as a column vector. 

The subset $L=(\overline\D\times\set{0})\cup\bigcup_{j=1}^m(J_j\times
\overline\D)$ of $\C^2$ is polynomially convex, so there is a smooth,
plurisubharmonic function $\rho_0$ on $\C^2$ such that $\rho_0=0$ on $L$
and $\rho_0>0$ on $\C^2\setminus L$.  Find $\delta>0$ small enough that
$\Phi_j(U_j''\times D_\delta^{n+3})$ is relatively compact in
$\Phi_0(D_r\times D_\epsilon^{n+3})$ for $j=1,\dots,m$.  By multiplying
$\rho_0$ by a large enough constant, we get $\rho_0\geq 1$ on each of
the sets $\partial D_r\times D_\epsilon$, $U_j''\times\partial D_s$, and
$\partial U_j''\times (D_s\setminus D_\delta)$.  The desired function
$\rho$ will be of the form $x\mapsto \rho_0(x_1,x_2)+\rho_1(x)$. 

We now define $\rho_1$ on the open subset
$$U=\Phi_0(D_r\times D_\epsilon^{n+3}) \cup \bigcup_{j=1}^m \Phi_j(U_j''
\times D_s^{n+3})$$
of $\C^4\times X$.  We start by choosing strictly positive, continuous, 
subharmonic functions $\alpha$, $\beta$ on $\C$, such that for
$j=1,\dots,m$, 
$$\alpha=\frac {\epsilon^2} 3 \text{ on } J_j, \qquad \alpha\geq 1
\text{ and } \beta\leq \frac 1 A \text{ near }\partial U_j'',$$
$$\beta\geq \frac 1 A - \eta \text{ on } U_j'\setminus U_j'',
\qquad \beta > \max\set{\alpha,A} \text{ on } 
D_r\setminus\bigcup_{j=1}^m U_j',
$$
where $\eta=\tfrac 1 A - 1 + \tfrac 1 3 (1-R)^2>0$. 
On each of the mutually disjoint open sets $\Phi_j(U_j''\times
D_s^{n+3})$, let
$$\rho_1(\Phi_j(x)) = \frac 1 3 |x_1|^2 + \frac 1 {\epsilon^2}
(\alpha(x_1) |x_2|^2 + |x''|^2),$$
where $x=(x_1,x')=(x_1,x_2,x'')$.  On the open set $\Phi_0((D_r\setminus
\bigcup_{j=1}^m \overline U_j')\times D_\epsilon^{n+3})$, let
$$\rho_1(\Phi_0(x))=\frac 1 3 |x_1|^2 + \frac 1
{\epsilon^2}\beta(x_1) |x'|^2.$$
By adding the set $\Phi_0((\bigcup_{j=1}^m \overline
U_j'\setminus U_j'') \times D_\epsilon^{n+3})$, we obtain a partition of
$U$, and on this set we let
$$\rho_1(\Phi_0(x))=\frac 1 3 |x_1|^2 + \frac 1 {\epsilon^2} \max\set{
\beta(x_1)|x'|^2, \alpha(x_1)|x_2|^2+|\pi_2 \Phi_j^{-1} 
\Phi_0(x)|^2},$$
where $\pi_2:\C^{n+4}\to\C^{n+2}$ chops off the first two coordinates. 
By the choice of $\alpha$ and $\beta$, we have
$$\beta(x_1)|x'|^2 \leq \alpha(x_1)|x_2|^2+|\pi_2 \Phi_j^{-1}
 \Phi_0(x)|^2 \qquad\text{for } x_1\text{ near } \partial U_j'',$$
and
$$\beta(x_1)|x'|^2 \geq \alpha(x_1)|x_2|^2+|\pi_2 \Phi_j^{-1}
 \Phi_0(x)|^2 \qquad\text{for } x_1 \text{ near }\partial U_j',$$
so $\rho_1$ is a well-defined, positive, continuous, strictly
plurisubharmonic function on $U$. 

On $K$, $\rho(x)=\tfrac 1 3 (|x_1|^2+|x_2|^2)$, so $\rho(K)=[0,\tfrac 2
3]$.  To complete the proof, we will verify that $\liminf\rho(x)\geq 1$
as $x$ goes to infinity in $U$, i.e., as $x$ eventually leaves each
compact subset of $U$.  This implies that $\rho$ exhausts the
neighbourhood $V=\rho^{-1}[0,1)$ of $K$. 

First let $\Phi_j(x)\to\infty_U$ with $x\in U_j''\times D_s^{n+3}$.  We
may assume that $x\notin U_j''\times D_\delta^{n+3}$ since the
$\Phi_j$-image of this set is relatively compact in $\Phi_0(D_r\times
D_\epsilon^{n+3})$.  Now either $x_1\to\partial U_j''$ with $x'\in
D_s^{n+3}\setminus D_\delta^{n+3}$, or $x_2\to\partial D_s$, or
$x''\to\partial D_s^{n+2}$, so $\liminf\rho(x)$ is larger than either
$\inf\rho_0(\partial U_j''\times(D_s\setminus D_\delta)\cup
U_j''\times\partial D_s)$ or $s^2/\epsilon^2$, so it is at least $1$.

Next let $\Phi_0(x)\to\infty_U$ with $x\in D_r\times D_\epsilon^{n+3}$,
$x_1\notin \bigcup_{j=1}^m U_j'$.  Then either $x_1\to\partial D_r$ or
$x'\to\partial D_\epsilon^{n+3}$, so $\liminf\rho(x)$ is larger than
either $\inf\rho_0(\partial D_r\times D_\epsilon)$ or 
$\inf\beta(D_r\setminus\bigcup_{j=1}^m U_j')$, so it is at least $1$.

Finally, let $\Phi_0(x)\to\infty_U$ with $x_1\in U_j'\setminus U_j''$. 
Then $x'\to\partial D_\epsilon^{n+3}$, and $\liminf\rho$ is at least
$\tfrac 1 3 (1 - \text{radius }U_j')^2 + \tfrac 1 A - \eta \geq 1$.
\qed \enddemo

Both hypotheses \therosteritem2 and \therosteritem3 may be reduced to
the affine case in the sense that they hold for a disc functional $H$ on
a complex manifold $X$ if they hold for every pullback $\phi^*H$ of $H$
by a holomorphic submersion $\phi$ from a domain of holomorphy in affine
space into $X$.  To see this for \therosteritem2, simply cover $X$ by
open balls and apply \therosteritem2 to the pullback of $H$ by the
inclusion of each ball into $X$.  The union of the covers thus obtained
for each ball is the required cover of $X$.  As for \therosteritem3,
take $h\in\A_X$, $w\in\T$, and
$\epsilon>0$.  Find $f\in\A_X$ such that $f(0)=h(w)$ and $H(f)\leq
EH(h(w))+\epsilon/2$.  Suppose $f$ and $h$ are both holomorphic on
$D_r$, $r>1$.  The union of the embedded discs $\set{(z,w,f(z))\,;\,z\in
D_r}$ and $\set{(0,z,h(z))\,;\, z\in D_r}$ is a Stein subvariety of
$D_r^2\times X$ and has a Stein neighbourhood $V$ by Siu's theorem.  Let
$\tau:V\to\C^m$ be an embedding and $\sigma:W\to\tau(V)$ be a
holomorphic retraction from a Stein neighbourhood $W$ of the submanifold
$\tau(V)$ in $\C^m$.  Now $\phi=\pr\circ\tau^{-1}\circ\sigma:W\to X$ is
a holomorphic submersion by which both $f$ and $h$ lift: $\tilde
f:z\mapsto\tau(z,w,f(z))$ is a lifting of $f$ and $\tilde
h:z\mapsto\tau(0,z,h(z))$ of $h$.  Here, $\pr:\C^2\times X\to X$ is the
projection.  To get condition \therosteritem3 for $H$ we postcompose by
$\phi$ the holomorphic maps provided by \therosteritem3 for $\phi^*H$
and the data $\tilde h$, $w$, $\epsilon/2$, and observe that
$$E\phi^*H(\tilde h(w))+\epsilon/2 \leq \phi^*H(\tilde f)+\epsilon/2 =
H(f)+\epsilon/2 \leq EH(h(w)) + \epsilon$$ 
by the choice of $f$. 

Suppose we have defined a class of disc functionals on each complex
manifold, preserved by taking pullbacks by holomorphic submersions. 
Consider the property of a disc functional of satisfying \therosteritem2
and \therosteritem3 and having a plurisubharmonic envelope.  The
Reduction Theorem implies that to establish this property for all disc
functionals of the given type, it suffices to prove it on domains of
holomorphy in affine space.  We do not know if the property of having a
plurisubharmonic envelope by itself reduces to the affine case.

\specialhead 2. Plurisubharmonicity of the Poisson and Riesz envelopes
\endspecialhead

\noindent Let $X$ be a complex manifold and $v:X\to\overline\R$ be a
Borel function which is locally bounded from above or from below.  Then
the Poisson functional $H_P^v$ on $X$ associated to $v$ is defined
by the formula 
$$H_P^v(f)=\int_\T v\circ f\,d\lambda, \qquad f\in\A_X.$$ 

Rosay's theorem \cite{Rs} states that the Poisson envelope $EH_P^v$ is
plurisubharmonic on $X$ when $v:X\to[-\infty,+\infty)$ is upper
semicontinuous, and then $EH_P^v$ is the largest plurisubharmonic
minorant of $v$.  This may also be derived from the Reduction Theorem
1.2 (to whose proof Rosay's ideas are essential) and plurisubharmonicity
of the Poisson envelope on domains in affine space, proved independently
(with different proofs) first by Poletsky \cite{P1} and later by Bu and
Schachermayer \cite{BS}. 
 
Now suppose $v:X\to(-\infty,+\infty]$ is plurisuperharmonic, not
identically $+\infty$.  Edigarian \cite{E1} has shown that the Poisson
envelope $EH_P^v$ is plurisubharmonic when $X$ belongs to a large class
of complex manifolds (the so-called class $\Cal P$ from \cite{LS1}), 
including domains in affine space.  Hence, $H_P^v$
satisfies hypothesis \therosteritem1 in the Reduction Theorem 1.2 (note
that the pullback $\phi^*H_P^v$ is the Poisson functional
$H_P^{v\circ\phi}$ of the plurisuperharmonic function $v\circ\phi$). 
Also, \therosteritem2 holds: take $U=X$ and $Z=v^{-1}(+\infty)$. 

Let us verify that $H_P^v$ satisfies \therosteritem3.  Take $h\in\A_X$,
$w\in\T$, and $\beta>EH_P^v(h(w))$.  By Lemma 9 in \cite{E1}, there is a
holomorphic map $G:D_r\times B_r\times V\to X$, where $B_r$ is the ball
of radius $r>1$ in $\C^n$, $n=\dim X$, and $V$ is a neighbourhood of
$h(w)$ in $X$, such that $G(0,y,x)=x$ for $y\in B_r$ and $x\in V$, and
for each $x\in V$, the average of $y\mapsto H_P^v(G(\cdot,y,x))$ on
$B_1$ is less than $\beta$.  If $J\subset U=h^{-1}(V)$ is a closed arc
in $\T$, then the average over $B_1\times J$ of $(y,w)\mapsto
H_P^v(G(\cdot,y,h(w)))$ is less than $\beta$.  Hence, there is $y_0\in
B_1$ (depending on $J$, of course) such that the average over $J$ of
$w\mapsto H_P^v(G(\cdot,y_0,h(w)))$ is less than $\beta$.  We let
$F:D_r\times U\to X$, $(z,w)\mapsto G(z,y_0,h(w))$, and observe that
$F(0,w)=G(0,y_0,h(w))=h(w)$ for $w\in U$. 

The Reduction Theorem 1.2 now implies the following result. 

\proclaim{2.1.  Theorem} Let $v:X\to(-\infty,+\infty]$ be a
plurisuperharmonic function, not identically $+\infty$, on a complex
manifold $X$.  Then the Poisson envelope $EH_P^v$ of $v$ is
plurisubharmonic on $X$ and is the largest plurisubharmonic minorant of
$v$.  \endproclaim

A plurisubharmonic function $u$ on a complex manifold $X$ defines the
Riesz functional $H_R^u$ on $X$ by the formula
$$H_R^u(f)=\frac 1 {2\pi}\int_\D\log|\cdot|\Delta(u\circ f),\qquad  
f\in\A_X,$$
where $\Delta(u\circ f)$ is considered as a positive Borel measure on
$\D$.  If $f\in\A_X$ and $u\circ f=-\infty$, then we set $H_R^u(f)=0$. 
By the Riesz representation formula,
$$EH_R^u = u + EH_P^{-u},$$
so we have the following immediate consequence of Theorem 2.1 and the
first part of Theorem 4.4 in \cite{LS1}.

\proclaim{2.2. Theorem}  The Riesz envelope $EH_R^u$ of a
plurisubharmonic function $u$ on a complex manifold $X$ is
plurisubharmonic on $X$ and is the largest negative plurisubharmonic
function on $X$ with Levi form no smaller than that of $u$.  
\endproclaim

We would like to draw the reader's attention to the important paper
\cite{P2}, in which Poletsky develops a new approach to the general
theory of disc functionals and proves, among other things, the
plurisubharmonicity of the Riesz envelope on domains in affine space.

\specialhead 3. Plurisubharmonicity of the Lelong envelope 
\endspecialhead

\noindent Let $\alpha$ be a positive function on a complex manifold
$X$.  We make no regularity assumptions on $\alpha$; it need
not even be measurable.  In \cite{LS1}, we defined the Lelong functional 
$H_L^\alpha$ on $X$ associated to $\alpha$ by the formula
$$H_L^\alpha(f)=\sum_{z\in\D}\alpha(f(z))\, m_z(f) \log|z|, \qquad
f\in\A_X.$$
By omitting the multiplicity $m_z(f)$, we obtain the {\it reduced Lelong
functional} $$\widetilde H_L^\alpha(f)=\sum_{z\in\D}\alpha(f(z))
\log|z|,$$ whose advantage over the Lelong functional is that
its pullback by a holomorphic map is a functional of the same type.  
(For this reason, and since the envelopes are the same, as proved below,
we suggest that the reduced Lelong functional be known as the Lelong
functional; the original Lelong functional is perhaps best forgotten.)

\proclaim{3.1.  Lemma} Let $\alpha$ be a positive function on a
domain $X$ in $\C^n$.  Then the Lelong functional $H_L^\alpha$ and
the reduced Lelong functional $\widetilde H_L^\alpha$ have the same
envelope.
\endproclaim

\demo{Proof} Since $H_L^\alpha\leq \widetilde H_L^\alpha$, it suffices to
show that for every $f\in\A_X$ and $\beta >H_L^\alpha(f)$, there exists
$g\in \A_X$ such that $g(0)=f(0)$ and $\widetilde H_L^\alpha(g)<\beta$.

By definition of $H_L^\alpha$, there exist finitely many distinct
points $a_1,\dots,a_N\in \D$ with multiplicities $m_j=m_{a_j}(f)$ such
that $\sum \alpha(f(a_j))m_j\log|a_j|<\beta$.  If this sum is $-\infty$, 
then there is $j$ with $a_j=0$ and $\alpha(f(a_j))>0$, so 
$\widetilde H_L^\alpha(f)=-\infty<\beta$.  Hence,
we may assume that $a_j\neq 0$ for all $j$. We set $m=m_1+\cdots+m_N$
and define $b=(b_0,\dots,b_m)\in \D^{m+1}$ as the vector
$(0,a_1,\dots,a_1,\dots,a_N,\dots,a_N)$ with the point
$a_j$ repeated $m_j$ times.  We will choose $t_j\in \T$ close to
$1$ such that $c_j=t_jb_j$ are all different, and show that
there exists $g\in \A_X$ such that $g(c_j)=f(b_j)$ for $j=0,\dots,m$.
Then
$$
\widetilde H_L^\alpha(g)\leq \sum_{j=1}^m \alpha(g(c_j))\log|c_j|
=\sum_{j=1}^m \alpha(f(b_j))\log|b_j|=
\sum_{j=1}^N \alpha(f(a_j))m_j\log|a_j|<\beta.
$$

We write $f=p+h$, where $p:\C\to \C^n$ is the polynomial map of degree 
at most $m$ that solves the interpolation problem $p(0)=f(0)$ and 
$p^{(k)}(a_j)=f^{(k)}(a_j)$ for $k=0,\dots,m_j-1$.  By the Newton 
interpolation formula,
$$
p(z)=f[b_0]+\sum_{j=1}^mf[b_0,\dots,b_j](z-b_0)\cdots(z-b_{j-1}) \tag{*}
$$ 
and $h(z)=f[b_0,\dots,b_m,z](z-b_0)\cdots(z-b_m)$, where 
$f[b_i,\dots,b_{i+j}]$ are the finite differences defined for $j=0$ as 
$f(b_i)$, and for $j>0$ as $f^{(j)}(b_i)/j!$ if $b_i=b_{i+j}$ and as
$\big(f[b_{i+1},\dots,b_{i+j}]-f[b_i,\dots,b_{i+j-1}]\big)
/\big(b_{i+j}-b_i \big)$ if $b_i\neq b_{i+j}$. 

Let $q:\C\to \C^n$ be the polynomial map of degree at most $m$ that 
solves the interpolation problem $q(c_j)=f(b_j)$ for $j=0,\dots,m$, 
let $k(z)=f[b_0,\dots,b_m,z](z-c_0)\cdots(z-c_m)$, and set $g=q+k$.
Then $g(c_j)=f(b_j)$ for $j=0,\dots,m$, so if we can show that
$g\to f$ uniformly on $\overline \D$ as $c\to b$, then $g\in \A_X$
for $c$ sufficiently close to $b$ and the proof is complete.  Clearly, 
$k\to h$ uniformly on $\overline \D$ as $c\to b$.  To prove that 
$q\to p$, we observe that the Newton formula for $q$ for interpolating 
the values of $g$ at the points $c_j$ is (*) with $p$, $f$, and $b$ 
replaced by $q$, $g$, and $c$, respectively.  We have $g[c_i]=f[b_i]$ 
for all $i$.  If $j>0$ and $b_i=b_{i+j}$, then clearly 
$g[c_i,\dots,c_{i+j}]=0$, and $f[b_i,\dots,b_{i+j}]=f^{(j)}(b_i)/j!=0$
since $b_i$ is one of the points $a_\nu$ and $j<m_\nu$.  It now follows 
easily by induction on $j$ that $g[c_i,\dots,c_{i+j}]\to 
f[b_i,\dots,b_{i+j}]$ for all $i$ and $j$, so $q\to p$ locally 
uniformly on $\C$, as $c\to b$. 
\qed
\enddemo

\proclaim{3.2.  Theorem} Let $\alpha$ be a positive function on a
complex manifold $X$.  Then the Lelong envelope $EH_L^\alpha$ and the
reduced Lelong envelope $E\widetilde H_L^\alpha$ associated to $\alpha$ are
equal and plurisubharmonic on $X$, and they are the largest negative
plurisubharmonic function on $X$ with Lelong numbers at least $\alpha$.
\endproclaim

\demo{Proof} If $X$ is a domain in affine space, then $EH_L^\alpha$ is
plurisubharmonic by Theorem 5.3 in \cite{LS1}, and equal to $E\widetilde
H_L^\alpha$ by Lemma 3.1.  In the general case, consider the class $\Cal
F_\alpha$ of negative plurisubharmonic functions on $X$ with Lelong
numbers at least $\alpha$.  It is easily seen that $\sup\Cal F_\alpha\leq
EH_L^\alpha\leq E\widetilde H_L^\alpha$.  By Proposition 5.1 in \cite{LS1},
which is stated for $H_L^\alpha$ but applies equally well to $\widetilde
H_L^\alpha$, the envelope $E\widetilde H_L^\alpha$ is plurisubharmonic if
and only if $E\widetilde H_L^\alpha\in\Cal F_\alpha$, and then $E\widetilde
H_L^\alpha=\sup\Cal F_\alpha$, so $E\widetilde H_L^\alpha=EH_L^\alpha$. 

It therefore remains to show that $E\widetilde H_L^\alpha$ is
plurisubharmonic on our complex manifold.  By the above, hypothesis
\therosteritem1 in the Reduction Theorem 1.2 is satisfied.  Hypothesis
\therosteritem2 is clear, since $\widetilde H_L^\alpha\leq 0$.  Finally,
by Lemma 2.3 in \cite{LS1}, $\widetilde H_L^\alpha$ (as well as
$H_L^\alpha$) satisfies the condition in Remark 1.3, so hypothesis
\therosteritem3 holds.  
\qed\enddemo

By Theorem 5.3 in \cite{LS1}, on a domain in a Stein manifold, the
Lelong envelope is the Poisson envelope of the so-called Lempert
function $k_\alpha$.  In fact, we deduced the plurisubharmonicity of the
Lelong envelope from this.  It is easy to see that this holds on all
manifolds, now that we know that $EH_P^{k_\alpha}$ and $EH_L^\alpha$ are
plurisubharmonic: $EH_P^{k_\alpha} \leq EH_L^\alpha$ \cite{LS1, p\. 
22}, and since $k_\alpha\geq EH_L^\alpha$ and $EH_L^\alpha$ is
plurisubharmonic, $EH_P^{k_\alpha} \geq EH_L^\alpha$.


\Refs

\widestnumber\key{MM}

\ref \key BS \by Bu, S. Q., and W. Schachermayer
\paper Approximation of Jensen measures by image measures under holomorphic
functions and applications 
\jour Trans. Amer. Math. Soc.  \vol 331  \yr 1992  \pages 585--608
\endref

\ref \key E1 \by Edigarian, A.
\paper A note on L\ac arusson-Sigurdsson's paper
\paperinfo preprint (2001), to appear in Math. Scand
\endref

\ref \key E2 \bysame
\paper A note on Rosay's paper
\paperinfo preprint  \yr 2001
\endref

\ref \key EP \bysame and E. A. Poletsky
\paper Product property of the relative extremal function 
\jour Bull. Polish Acad. Sci. Math. \vol 45 \yr 1997 \pages 331--335
\endref

\ref \key LS1 \by L\ac arusson, F., and  R. Sigurdsson
\paper Plurisubharmonic functions and analytic discs on manifolds
\jour J. reine angew. Math.  \vol 501  \yr 1998  \pages 1--39
\endref

\ref \key LS2 \bysame
\paper Plurisubharmonic extremal functions, Lelong numbers and coherent ideal
sheaves 
\jour Indiana Univ. Math. J.   \vol 48  \yr 1999  \pages 1513--1534
\endref

\ref \key P1 \by Poletsky, E. A.
\paper Plurisubharmonic functions as solutions of variational problems 
\inbook Several complex variables and complex geometry (Santa Cruz, CA, 1989)
\pages 163--171  \bookinfo Proc. Sympos. Pure Math., 52, Part 1 
\publ Amer. Math. Soc. \yr 1991
\endref

\ref \key P2 \bysame
\paper The minimum principle
\paperinfo preprint (1999), to appear in Indiana Univ. Math. J 
\endref

\ref \key PS  \bysame and B. V. Shabat
\paper Invariant metrics
\inbook Several Complex Variables III  \bookinfo Encyclopaedia of
Mathematical Sciences, volume 9  \publ Springer-Verlag \yr 1989
\pages 63--111
\endref

\ref \key Rs \by Rosay, J.-P.
\paper Poletsky theory of disks on holomorphic manifolds
\paperinfo preprint (2001), to appear in Indiana Univ. Math. J 
\endref

\ref \key Ry  \by Royden, H. L.
\paper The extension of regular holomorphic maps
\jour Proc. Amer. Math. Soc.  \vol 43  \yr 1974  \pages 306--310
\endref

\endRefs
\enddocument